\documentclass[12pt]{article}
\newcommand{\copyleft}{
GNU FDL\thanks{
Copyright (C) 2009 Peter G. Doyle.
Permission is granted to copy, distribute and/or modify this document
under the terms of the GNU Free Documentation License,
as published by the Free Software Foundation;
with no Invariant Sections, no Front-Cover Texts, and no Back-Cover Texts.
}}
\title{The Kemeny constant of a Markov chain}
\author{Peter Doyle}
\date{Version 1.0 dated 14 September 2009
\\ \copyleft}
\usepackage{url}
\usepackage{tensor}

\newcommand{\comment}[1]{}

\newcommand{\inv}{{-1}}

\newcommand{\Pinf}{P^{\infty}}
\newcommand{\Prob}{\mathrm{Prob}}
\newcommand{\Exp}{\mathrm{Exp}}
\newcommand{\Var}{\mathrm{Var}}

\newcommand{\eps}{\epsilon}
\newcommand{\tr}{\mathrm{Tr}}
\newcommand{\qed}{\rule{2mm}{2.5mm}}

\newcommand{\Pij}{\tensor{P}{_i^j}}
\newcommand{\Zij}{\tensor{Z}{_i^j}}
\newcommand{\Zjj}{\tensor{Z}{_j^j}}
\newtheorem{theorem}{Theorem}
\newtheorem{prop}[theorem]{Proposition}
\newtheorem{corollary}[theorem]{Corollary}

\begin{document}

\maketitle

\begin{abstract}
Given an ergodic finite-state Markov chain,
let $M_{iw}$ denote the
mean time from $i$ to equilibrium,
meaning the expected time, starting from $i$,
to arrive at a state selected
randomly according to the equilibrium measure $w$ of the chain.
John Kemeny observed that
$M_{iw}$ does not depend on starting the point $i$.
The common value $K=M_{iw}$ is
the \emph{Kemeny constant} or \emph{seek time}
of the chain.
$K$ is a spectral invariant, to wit,
the trace of the resolvent matrix.
We review basic facts about the seek time,
and connect it to the bus
paradox and the Central Limit Theorem for ergodic Markov chains.
\end{abstract}

\centerline{\emph{For J. Laurie Snell}}
\subsection*{The seek time}

We begin by reviewing basic facts and establishing notation for Markov chains.
For background, see
Kemeny and Snell
\cite{kemenysnellfinite} or
Grinstead and Snell
\cite{Grinsteadsnell},
bearing in mind that the notation here is somewhat different.

Let $P$ be the transition matrix of an ergodic finite-state Markov chain.
We write the entries of $P$ using tensor notation, with
$\Pij$ being the probability that from state $i$ we move to state $j$.
(There is some possibility here of confusing superscripted indices with
exponents, but in practice it should be clear from context which is meant.)
The sequence of matrix powers $I,P,P^2,P^3,\ldots$ has a (Cesaro) limit which we will
denote by $\Pinf$.
We have $P \Pinf = \Pinf P = \Pinf$.
The rows of $\Pinf$ are identical:
\[
\tensor{(\Pinf)}{_i^j} = w^j
.
\]
Like the rows of $P$, the row vector $w$ is a probability distribution:
$w^j \geq 0$, $\sum_j w^j = 1$.
$w$ is the \emph{equilibrium distribution} of the chain.
The entry $w^j$ tells the steady-state probability that the chain is in
state $j$.

The row vector $w$ is, up to multiplication by a constant,
the unique row vector fixed by $P$:
\[
\sum_i w^i \Pij = w^j
.
\]
That makes it a row eigenvector corresponding to the eigenvalue $1$.
The corresponding column eigenvector is the constant vector:
\[
\sum_j \Pij  \cdot 1 = 1
.
\]

Let $M_{ij}$ be the expected time to get from state $i$ to state $j$,
where we take $M_{ii}=0$.
The \emph{mean time from $i$ to equilibrium} is
\[
M_{iw} = \sum_j M_{ij} w^j
.
\]
This tells the expected time to get from state $i$ to a state $j$ selected
randomly according to the equilibrium measure $w$.

John Kemeny
(see \cite[4.4.10]{kemenysnellfinite}, \cite[p. 469]{Grinsteadsnell})
observed:
\begin{theorem}
$M_{iw}$ doesn't depend on $i$.
\end{theorem}

The common value of the $M_{iw}$'s, denoted $K$,
is the \emph{Kemeny constant} or \emph{seek time}
of the chain.

{\bf Proof.}
Observe that the function $M_{iw}$ is
\emph{discrete harmonic}, meaning that it has the \emph{averaging property}
\[
\sum_j \Pij M_{jw} =M_{iw}
.
\]
The reason is that
taking a step away from $i$ brings you one step closer to your destination,
except when your destination is $i$ and the step begins a wasted journey from
$i$ back to $i$:
This happens with probability $w^i$,
and the expected duration of the wasted journey is
$\frac{1}{w^i}$,
because the mean time between visits to $i$ is the reciprocal of the
equilibrium probability of being there.
Thus
\[
M_{iw} -1 + w^i \frac{1}{w^i} = \sum_j \Pij  M_{jw}
,
\]
so
\[
M_{iw} = \sum_j \Pij M_{jw}
.
\]

But now by the familiar \emph{maximum principle},
any function $f_i$ satisfying
\[
\sum_j \Pij f_j = f_i
\]
must be constant:
Choose $i$ to maximize $f_i$, and observe that the maximum
must be attained
also for any $j$ where $\Pij>0$;
push the max around until it is attained everywhere.
So $M_{iw}$ doesn't depend on $i$. $\qed$

{\bf Note.}
The application of the maximum principle we've made here
shows that the only column eigenvectors having eigenvalue $1$
for the matrix $P$
are the constant vectors---a fact that was stated not quite explicitly above.

The foregoing argument shows the mean time from $i$ to equilibrium is
constant---but what is its value?
For this we return to
Kemeny's original proof of constancy for $M_{iw}$,
which involved writing an explicit formula
for $M_{iw}$, and noticing that it doesn't depend on $i$.

Define the \emph{resolvent} or \emph{fundamental matrix} or \emph{Green's function}
\begin{eqnarray*}
Z
&=&
(I-\Pinf) + (P- \Pinf) + (P^2-\Pinf) + \ldots
\\&=&
(I-(P-\Pinf))^\inv - \Pinf
.
\end{eqnarray*}
Please be aware that this resolvent $Z$
differs from the variant
used by Kemeny and Snell
\cite{kemenysnellfinite}
and
Grinstead and Snell
\cite{Grinsteadsnell},
which with our notation would be 
$(I-(P-\Pinf))^\inv$.
As others have observed
(cf. Meyer \cite{Meyer}; Aldous and Fill \cite{Aldousfill}),
for the version of $Z$ we use here,
the entries $\Zij$ have the natural
probabilistic interpretation as the `expected excess visits to $j$,
starting from $i$, as compared with a chain started in equilibrium'.
Accordingly we have
\[
\sum_i w^i \Zij = 0
\]
and
\[
\sum_j \Zij = 0
.
\]

Since $\Zij$ measures excess visits to $j$
starting at $i$, relative to starting in equilibrium,
we obviously have $\Zjj \geq \Zij$,
because to make excess visits to $j$ starting
from $i$ you first have to get to $j$.
And the discrepancy $\Zjj-\Zij$ is just
$M_{ij} w^j$, because in equilibrium this is the expected number of visits
to $j$ over an interval of expected length $M_{ij}$.
From this we get the familiar formula
\[
M_{ij} = 
(\Zjj - \Zij)
\frac{1}{w^j}
.
\]
(Cf. \cite[4.4.7]{kemenysnellfinite}, \cite[p. 459]{Grinsteadsnell})

\begin{prop}[\protect{Kemeny and Snell \cite[4.4.10]{kemenysnellfinite}}]
Kemeny's constant is the trace of the resolvent $Z$:
\[
K = M_{iw} = \sum_j \Zjj
.
\]
\end{prop}

{\bf Proof.}
\begin{eqnarray*}
M_{iw}
&=&
\sum_j M_{ij} w^j
\\&=&
\sum_j
\Zjj
-
\sum_j
\Zij
\\&=&
\sum_j
\Zjj
,
\end{eqnarray*}
using the fact that
$\sum_j \Zij = 0$.
$\quad \qed$

This formula provides a
computational verification that Kemeny's constant is constant,
but doesn't explain \emph{why} it is constant.
Kemeny felt this keenly:  A prize was offered for a more `conceptual'
proof, and awarded---rightly or wrongly---on the basis of
the maximum principle
argument outlined above.

Still, there are advantages to having an explicit formula.
For starters, the explicit formula reveals that
the seek time is a spectral
invariant of the matrix $I-P$:
If we denote the eigenvalues of $I-P$ by
$\lambda_0=0,\lambda_1,\ldots,\lambda_{n-1}$,
then the eigenvalues of $Z= (I-(P-\Pinf))^\inv-\Pinf$ are
$0,\frac{1}{\lambda_1},\ldots,\frac{1}{\lambda_{n-1}}$,
and
\[
K=\tr(Z) = \frac{1}{\lambda_1} + \ldots + \frac{1}{\lambda_{n-1}}
.
\]
(Cf.
Meyer \cite{Meyer},
Aldous and Fill \cite{Aldousfill},
Levene and Loizou \cite{Leveneloizou}.)
In terms of the eigenvalues
$\alpha_0=1,\alpha_1=1-\lambda_1,\ldots,
\alpha_{n-1}=1-\lambda_{n-1}$
of $P$ we have
\[
K=\tr{Z} = \frac{1}{1-\alpha_1}+\ldots+\frac{1}{1-\alpha_{n-1}}
.
\]
We'll have more to say about this later.

\section*{Kemeny on the bus}

We now explore connections between the seek time, the famous
\emph{bus paradox},
and the Central Limit Theorem (CLT) for Markov Chains.

Just to recall, the bus paradox is that for a Poisson
process like radioactive decay,
the expected length of the interval between
the events that bracket any given instant
is twice the average interval between
events---so if buses are dispatched by a Geiger counter, you must expect
to wait twice as long for the bus as you would if the buses came at
regular intervals.
The explanation for this is 
that any given instant is more likely to land in a long inter-event
interval than in a
short inter-event interval,
so we're taking a weighted average of the intervals between
events, emphasizing the longer intervals,
and this makes the expected waiting time longer than the average
inter-event time.
This inequality will be true for any renewal process in equilibrium:
The factor of $2$ disparity is special
to the Poisson process, and arises because
the Poisson process is memoryless and time-reversible.

To make the connection of the seek time
to the bus paradox, we think about the
\emph{mean time from equilibrium to $j$}:
\[
M_{wj} = \sum_i w^i M_{ij}
.
\]
Here we choose a starting state $i$ at random according to $w$,
and see how long it takes to get to $j$.
This is backwards from what we did to
define the seek time $M_{iw}$,
where we looked at the time to get from $i$ to equilibrium.
%$M_{wj}$ has been known to go by the name of the \emph{preKemeny nonconstant}.

Unlike $M_{iw}$, which is independent of $i$,
the quantity $M_{wj}$ depends on the state $j$.
Choosing the target state $j$ at random according to $w$ gets
us back to $K$:
\[
\sum_j M_{wj} w^j = K
.
\]

\begin{prop}[\protect{Kemeny and Snell \cite[4.4.9]{kemenysnellfinite}}]
The mean time from equilibrium to $j$ is
\[
M_{wj} = \Zjj \frac{1}{w^j}
.
\]
\end{prop}

{\bf Proof.}
\[
M_{wj}
= 
\sum_i w^i M_{ij}
= 
\sum_i w^i (\Zjj-\Zij)\frac{1}{w^j}
=
\Zjj\frac{1}{w^j}
,
\]
because $\sum_i w^i =1$ and $\sum_i w^i \Zij = 0$.
$\quad \qed$

To take advantage of this formula for $M_{wj}$,
and specifically,
to use it to derive the CLT for Markov chains,
we now recall a bit of renewal theory.
(Cf. Feller \cite[Chapter XIII]{Feller1}.)

A discrete renewal process is effectively just what you get if you
watch a discrete-time Markov chain (possibly having infinitely many states)
and take note of the times at which it is in some
fixed state $a$.
These times are called \emph{renewal times} or \emph{renewals}.
The name derives from the fact
that each time the chain reaches $a$ it begins anew.
The variables that tell the elapsed time between successive renewals
are independent and identically distributed.

So let $X$ be a random variable whose distribution is that for the time
between successive renewals,
and let it have mean $\mu$ and variance $\sigma$.
We want to express the mean time $\tau$ from equilibrium to the next
renewal in terms of $\mu$ and $\sigma$.

{\bf Proposition (The bus equality).}
For a discrete renewal process with interarrival times having mean $\mu$ and
variance $\sigma^2$,
the mean time $\tau$ from equilibrium to the next
renewal is
\[
\tau = \frac{\mu^2+\sigma^2}{2\mu} - \frac{1}{2}
.
\]

The term $-\frac{1}{2}$ here is an artifact of using discrete time.

{\bf Proof.}
Let
\[
p(n) = \Prob(X=n)
,
\]
so
\[
\mu = \Exp(X) = \sum_n n p(n)
\]
and
\[
\mu^2 + \sigma^2 = \Exp(X^2) = \sum_n n^2 p(n)
.
\]
The expected time $\tau$ from equilibrium to the next renewal is
\begin{eqnarray*}
\tau
&=&
\frac{\sum_n \frac{n-1}{2} n p(n)}{\sum_n n p(n)}
\\&=&
\frac{1}{2} \left (
\frac{\sum_n n^2 p(n)}{\sum_n n p(n)} - 1
\right)
\\&=&
\frac{\mu^2+\sigma^2}{2\mu} - \frac{1}{2}
.
\quad \qed
\end{eqnarray*}

\begin{corollary}[The bus inequality]
\[
\tau \geq \frac{\mu}{2} - \frac{1}{2}
,
\]
with equality just if $\sigma=0$ (all interarrival times equal).
$\quad \qed$
\end{corollary}

As in the bus equality above,
the term $-\frac{1}{2}$ here is an artifact of using discrete time.

The bus equality shows that knowing the time $\tau$
from equilibrium to the next renewal
is equivalent to knowing
$\sigma^2$, the variance of the renewal time:
\[
\tau = \frac{\mu^2+\sigma^2}{2\mu} - \frac{1}{2}
;
\]
\[
\sigma^2 = 2 \mu \tau + \mu - \mu^2
.
\]
Of course the mean renewal time $\mu$ is involved here, too:  We take
that for granted.

Now let's return to our Markov chain,
and take for our renewal process visits to a given state $j$.
The mean time between renewals is $\mu=\frac{1}{w^j}$.
The expected time in equilibrium to the next renewal is
$\tau=M_{wj} = \Zjj\frac{1}{w^j}$.
But above we saw that
\[
\sigma^2 = 2 \mu \tau + \mu - \mu^2
\]
so
\[
\sigma^2 = 2 \Zjj\frac{1}{(w^j)^2} + \frac{1}{w^j} - \frac{1}{(w^j)^2}
.
\]
Going back the other way, from $\mu$ and $\sigma^2$ we can find $\Zjj$
(cf. Feller \cite[(5.1) on p. 443]{Feller1}):
\[
\Zjj = \frac{\sigma^2-\mu+\mu^2}{2 \mu^2}
.
\]

Another piece of information about a renewal process that
is equivalent to knowing $\sigma^2$ or $\tau$
(or $\Zjj$, in the Markov chain case we just discussed)
is the variance for the number of
renewals over a long period, which shows up in the Central Limit Theorem
for renewal processes:

\begin{theorem}[CLT for renewal processes]
For a renewal process whose renewal time has mean $\mu$ and variance
$\sigma^2$,
the number of renewals over a long time $T$ is approximately Gaussian
with mean $T \frac{1}{\mu}$ and variance $T \frac{\sigma^2}{\mu^3}$.
\end{theorem}

{\bf Idea of proof.}
To see a large number $N$ of renewals will take time
\[
T \approx N \mu \pm \sqrt{N} \sigma
.
\]
so the density of renewals over this interval is
\begin{eqnarray*}
\frac{N}{T}
&=&
\frac{N}{N \mu \pm \sqrt{N} \sigma}
\\&=&
\frac{1}{\mu \pm \frac{1}{\sqrt{N}}\sigma}
\\&=&
\frac{1}{\mu(1\pm\frac{1}{\sqrt{N}}\frac{\sigma}{\mu})}
\\&\approx&
\frac{1}{\mu} (1 \pm \frac{1}{\sqrt{N}}\frac{\sigma}{\mu})
\\&=&
\frac{1}{\mu} \pm \frac{1}{\sqrt{N}}\frac{\sigma}{\mu^2}
\\&=&
\frac{1}{\mu} \pm \frac{1}{\sqrt{N \mu}}\frac{\sigma}{\mu^{\frac{3}{2}}}
\\&\approx&
\frac{1}{\mu} \pm \frac{1}{\sqrt{T}}\frac{\sigma}{\mu^{\frac{3}{2}}}
\\&=&
\frac{1}{\mu} \pm \frac{1}{\sqrt{T}}\sqrt{\frac{\sigma^2}{\mu^3}}
.
\end{eqnarray*}
Thus
\[
N \approx T \frac{1}{\mu} \pm \sqrt{T \frac{\sigma^2}{\mu^3}}
. \quad \qed
\]

{\bf Note.}
Feller
\cite[p. 341]{Feller1}
gives the following formulas
for $\Exp(N)$ and $\Exp(N^2)$,
which he derives---or rather, asks readers to derive---using
generating functions:
\[
\Exp(N) = 
\frac{T+1}{\mu} + \frac{\sigma^2-\mu-\mu^2}{2 \mu^2} + o(1)
.
\]
\[
\Exp(N^2) =
\frac{(T+2)(T+1)}{\mu^2} + \frac{2 \sigma^2-2\mu-\mu^2}{\mu^3}T + o(T)
.
\]
These combine to give
\[
\Var(N) = \Exp(N^2) - Exp(N)^2 = \frac{\sigma^2}{\mu^3} T + o(T)
,
\]
which is the same as what we get from the CLT.
%Except that we could replace o(T) with O(1), which we leave it to readers
%to notice for themselves.

The CLT for renewal processes translates into the
following special case of the CLT for Markov chains
(special because we are considering only the number of visits to one
particular state, not the long-term average of a general function of
the state).
\begin{corollary}
For an ergodic Markov chain
with resolvent $Z$, the number of visits
to $j$ over a long time is approximately Gaussian with mean
$T \frac{1}{w^j}$ 
and variance
\[
T \frac{\sigma^2}{\mu^3} = T(2 \Zjj w^j + (w^j)^2 -w^j)
. \quad \qed
\]
\end{corollary}

Grinstead and Snell
\cite[p. 466]{Grinsteadsnell}
attribute this formula for the variance to Frechet.

Now just as in the case of the bus inequality, we get information
from the fact that the variance here must be $\geq 0$:
\[
2 \Zjj w^j + (w^j)^2 -w^j \geq 0
,
\]
so
\[
\Zjj \geq \frac{1 - w^j}{2}
.
\]
Summing over $j$ gives an inequality for the seek time $K$:
\begin{prop}
\[
K = \sum_j \Zjj \geq \frac{n-1}{2}
.
\quad \qed
\]
\end{prop}

This inequality for the seek time
was observed by 
Levene and Loizou
\cite{Leveneloizou}.
They derived it from the fact that if the non-1 eigenvalues of $P$ are
$\alpha_k$, $1 \leq k \leq n-1$
then the non-1 eigenvalues of $I-P+P^\infty$
are $1-\alpha_k$
and the non-0 eigenvalues of $Z = (I-P+P^\infty)^\inv-P^\infty$
are $\lambda_k=\frac{1}{1-\alpha_k}$.
But the $\alpha_k$'s lie in the unit disk,
which maps to the region $\{x+iy : x \geq \frac{1}{2}\}$ under the map
taking $z$ to $\frac{1}{1-z}$,
so the non-0 eigenvalues of $Z$ have real part $\geq \frac{1}{2}$,
and thus $K$, which is real, satisfies
\[
K = \tr(Z) =
\sum_{k=1}^{n-1} \lambda_k =
\sum_{k=1}^{n-1} \frac{1}{1-\alpha_k} \geq \frac{n-1}{2}
.
\]

{\bf Taking stock.}
From the resolvent $Z$ we've computed the variance of 
the \emph{return time}, meaning the time get from a designated
starting state $j$ back to $j$.
If instead we're interested in the variance of a
\emph{hitting time}, meaning the time to get from $i$ to $j \neq i$,
we'll need to look at $Z^2$.
We'd need this for the general CLT for Markov chains,
which as noted above deals with the long-term average of a general function
of the state of the chain, and requires knowing
the covariance of the number of visits to a pair of states $i$ and $j$.
Looking beyond mean and variance,
we get $k$th moments of return times from $Z^{k-1}$
and $k$th moments of hitting times from $Z^k$.
This was already in evidence for first moments:
We need $Z$ to find mean hitting times $M_{ij}=(\Zjj-\Zij)\frac{1}{w^j}$,
whereas
expected return times $\frac{1}{w^j}$ don't require knowing $Z$ at all.

% end comment

\bibliographystyle{hplain}
\bibliography{kc}

\end{document}